# Error-Correcting Codes Derived from Combinatorial Games

*Aviezri S. Fraenkel*

Department of Applied Mathematics & Computer Science
The Weizmann Institute of Science
Rehovot 76100, Israel
Email:   fraenkel@wisdom.weizmann.ac.il
WWW HomePage:   http://www.wisdom.weizmann.ac.il/~fraenkel/fraenkel.html

ABSTRACT.   The "losing positions" of certain combinatorial games constitute linear error detecting and correcting codes. We show that a large class of games that can be cast in the form of *annihilation games*, provides a potentially polynomial method for computing codes (*anncodes*). We also give a short proof of the basic properties of the previously known *lexicodes*, which are defined by means of an exponential algorithm, and are related to game theory. The set of lexicodes is seen to constitute a subset of the set of anncodes. In the final section we indicate, by means of an example, how the method of producing lexicodes can be applied optimally to find anncodes. Some extensions are indicated.

**1. Introduction.**   Connections between combinatorial games (simply *games* in the sequel) and linear error-correcting codes (*codes* in the sequel) have been established by Conway and Sloane [1986], Conway [1990], Brualdi and Pless [1993], where lexicodes, and some of their connections to games, are explored. *Lulei d'mist'fina*, we would venture to say that perhaps lexicodes *per se* are not the best (at least not the only) tool to appreciate this connection and to exploit it, since the construction of lexicodes involves exponential computations. Our aim is to extend the connection between games and codes to a large class of games, and to formulate a potentially polynomial method for generating codes from games. We also establish the basic properties of lexicodes by a simple, transparent method.

Let $\Gamma$, any finite digraph, be the *groundgraph* on which we play the following general 2-player game. Initially, distribute a positive finite number of tokens on the vertices of $\Gamma$. Multiple occupation is permitted. A move consists of selecting an occupied vertex and moving a single token from it to a neighboring vertex along a directed edge, occupied or not. The player first unable to move loses and the opponent wins. If there is no last move, the play is declared a draw. It is easy to see (since $\Gamma$ is finite), that a draw can arise only if $\Gamma$ is *cyclic*, i.e. , $\Gamma$ has cycles or loops. Games in this class—which includes Nim and Nim-like games for the case where $\Gamma$ is acyclic—have polynomial strategies, in general (Fraenkel [$\geq$ 1997]). It turns out that the $P$-positions (defined below) of any game in this class constitute a code.

It further turns out that, if $\Gamma$ is cyclic, then the structure of the $P$-positions is much richer if the above described game is replaced by an *annihilation game* (*anngame* for short). When a token is moved onto a vertex $u$ in this game, the number of tokens on $u$ is reduced mod 2, i.e., if the number of tokens on $u$ is even (odd) after the move, we remove all (remove all but one) tokens from $u$.

If $\Gamma$ is acyclic, it is easy to see by game-strategy considerations (or using the $g$-function defined below), that the strategies of the non-annihilation game and the anngame are identical, so both have the same $P$-positions—at most the length of play is affected. Thus, for the sake of a unified treatment, we may as well assume that all our games are anngames.

Summarizing, we assume, without loss of generality, that our above defined class of games consists of anngames. That is, given a finite digraph $\Gamma$, the groundgraph. Initially distribute tokens on the vertices, at most one token per vertex. A move consists of selecting an occupied vertex and moving its token to a neighboring vertex $u$ along a directed edge. If $u$ was occupied prior to this move, then the incoming and resident tokens on $u$ are both annihilated, i.e., they both disappear from play. The player first unable to move loses and the opponent wins. If there is no last move, the outcome is a draw for both players.

With a groundgraph $\Gamma$ on which an anngame $A$ is played, we associate its *annihilation graph* $G = (V, E)$, called *anngraph* for short, in which $V$ is the set of positions of $A$; and if $u, v \in V$ then $(u, v) \in E$ if and only if there is a move from $u$ to $v$ in $A$. The basic facts are that $G$ is a vector space over $GF(2)$. The *generalized Sprague-Grundy function* $\gamma$ on $G$, is a homomorphism from $V^f$, the linear subspace of $V$ on which $\gamma$ is finite, onto $GF(2)^t$—the set of binary vectors of size $t$, which is identified with the set of numbers $\{0, 1, \ldots, 2^t - 1\}$—for some $t \in \mathbb{Z}^0$. The kernel $V_0$, the subspace on which $\gamma = 0$, is also the set of $P$-positions of the annihilation game. The quotient space $V/V_0 = \{V_i : 0 \leq i < 2^t\}$ is the coset $V_i$ on which $\gamma = i$. This gives very precise information about the structure of $G$: its maximum finite $\gamma$-value is a power of 2 less 1, and the sets on which $\gamma = i$ have the same size for all $i \in \{0, \ldots, 2^t - 1\}$. Moreover, $V_0$ constitutes an *anncode* (annihilation game code).

If $\Gamma$ is cyclic, then generally $\gamma \neq g$ and $A$ has a distinctly different character and strategy than in the case when $\Gamma$ is acyclic. Though $G$ has $2^n$ vertices, all the relevant information can be extracted from an induced subgraph of size $O(n^4)$, by an $O(n^6)$ algorithm, which is often much more efficient. Annihilation games were suggested by John Conway. The above stated facts can be found in Fraenkel [1974], Fraenkel and Yesha [1976, 1979, 1982] (especially the latter), Yesha [1978], and Fraenkel, Tassa and Yesha [1978]. Ferguson [1984] considered misère annihilation play, i.e., the player first unable to move wins, and the opponent loses. A more transparent presentation of annihilation games is to appear in a forthcoming book, by Fraenkel ($\geq$1997), where also $\gamma$ is defined and illustrated. See also Fraenkel and Yesha [1986] for $\gamma$, which was defined by Smith [1966], and expounded in Fraenkel and Perl [1975].

In §2 we bring a number of examples, designed to illustrate connections between games, anncodes and lexicodes; exponential and polynomial digraphs and computations associated with them. In §3 we give a short proof that $g$ is linear on the lexigraph associated with lexicodes, leading to the same kind of homomorphism that exists for anncodes. Some natural further questions are posed at the end of §3, including the definition of anncodes over $GF(q)$, $q \geq 2$. In §4 we indicate, by means of a larger example, how a greedy algorithm applied to an anncode can reduce a computation of a code by a factor of 4,000 compared to a similarly computed lexicode.

The anncode method is potentially polynomial whereas the lexicode method is exponential. But it is too early yet to say to what extent the potential of the anncode method can be realized for producing new efficient codes.

**2. Examples.** Given a finite digraph $G = (V, E)$, we define, for any $u \in V$, the set of *followers* $F(u)$ and *ancestors* $F^{-1}(u)$ by

$$F(u) = \{v \in V : (u,v) \in E\}, \qquad F^{-1}(u) = \{w \in V : (w,u) \in E\}.$$

If the vertices of $G$ are game positions and the edges moves, we define, as usual, a *P-position* (*N-position*) of the game as any position from which the *P*revious (*N*ext) player can win, no matter how the opponent plays, subject to the rules of the game. Denote by $\mathcal{P}$ ($\mathcal{N}$) the set of all *P*-positions (*N*-positions) of a game. The following basic relationships hold:

$$u \in \mathcal{P} \quad \text{iff} \quad F(u) \subseteq \mathcal{N}, \qquad u \in \mathcal{N} \quad \text{iff} \quad F(u) \cap \mathcal{P} \neq \emptyset.$$

(If $G$ has cycles or loops, then the game may also contain (dynamic) $D$ (*D*raw)-positions, which satisfy:

$$u \in \mathcal{D} \quad \text{iff} \quad F(u) \subseteq \mathcal{D} \cup \mathcal{N}, \qquad F(u) \cap \mathcal{D} \neq \emptyset.)$$

To understand the examples below we don't need $\gamma$ or $g$; it suffices to note that $\mathcal{P}$ is the set of vertices on which $\gamma$ ($g$) is 0, and $\mathcal{P}$ can be recognized by purely game-theoretic considerations, as the set on which the previous (second) player can win. In all these examples, we play an annihilation game $A$ on the given groundgraphs $\Gamma$.

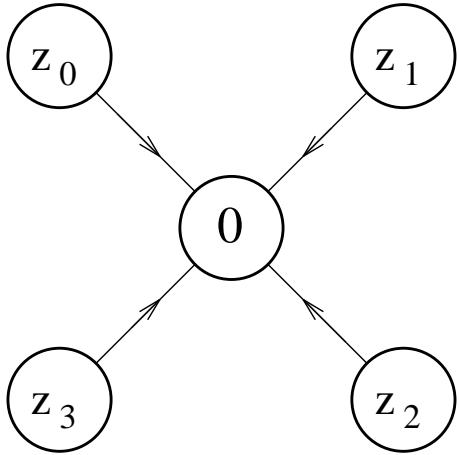

Fig. 1. An acyclic groundgraph for annihilation.

EXAMPLE 1. Let $\Gamma$ be the digraph depicted in Fig. 1. It is easy to see that with an odd number of tokens on the $z_i$ player I can win, and with an even number player II can win in $A$ played on $\Gamma$.

In this and in the following examples, think of the $z_i$ as unit vectors of a vector space $V$ of dimension $n$, where $n-1$ is the largest index of the $z_i$ (Fraenkel and Yesha [1982]). In the present example, $z_0 = (0001), \ldots, z_3 = (1000)$. Encoded by the unit vectors, we then have,

$$\mathcal{P} = \{(0000), (0011), (0101), (0110), (1001), (1010), (1100), (1111)\},$$

or, encoded in decimal, our anncode is given by $\mathcal{P} = \{0, 3, 5, 6, 9, 10, 12, 15\}$.

Note that $\mathcal{P}$ is a linear code with *minimal distance* $d = 2$.

EXAMPLE 2. Consider $A$ played on $\Gamma$ given in Fig. 2. If $z_0$ and $z_1$ host a token each, then any move causes annihilation. Therefore the position consisting of one token each on $z_0$, $z_1$, $z_2$ (or $z_3$ instead of $z_2$), is a $P$-position. Using our decimal encoding, we then see that $\mathcal{P} = \{0, 7, 11, 12\}$, which is also a linear code with $d = 2$.

EXAMPLE 3. Consider $A$ played on a Nim-heap of size 5, i.e., $\Gamma$ consists of the leaf 0 and the vertices $z_0, \ldots, z_4$, where $(z_j, z_i) \in E(\Gamma)$ if and only if $i < j$. It is not hard to see that then $\mathcal{P} = \{0, 7, 25, 30\}$, which is an anncode with $d = 3$. Note that precisely the same code is given by the $P$-positions of the annihilation game $A$ played on the ground graph $\Gamma$ of Fig. 3.

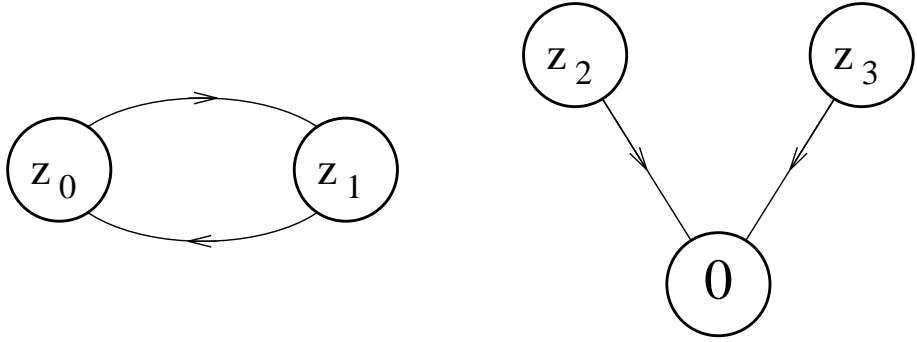

FIG. 2. A cyclic groundgraph.

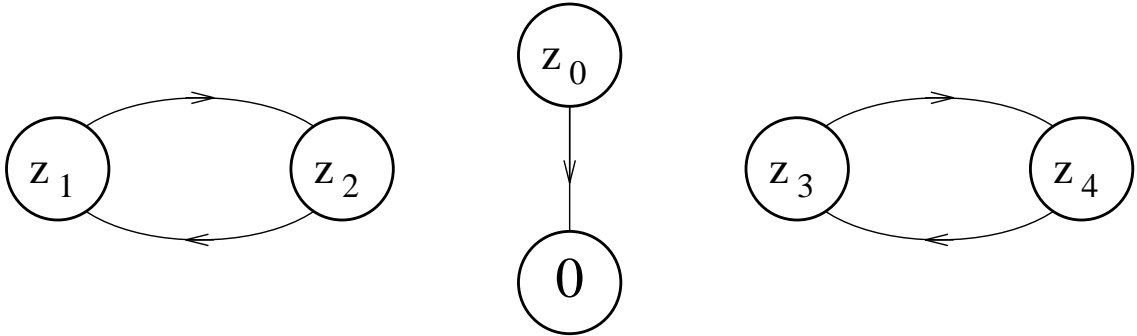

FIG. 3. Another cyclic groundgraph.

To continue with our examples, we now define lexicodes precisely. This is also needed for §3.

Let $W$ be an $n \times n$ matrix over $GF(2)$, of rank $\geq m$, whose columns $1, \ldots, m$, counted from the right, constitute a basis of $V^m$, the $m$-dimensional vector subspace of $V^n$ ($m \leq n$) over $GF(2)$. Then there are rows $1 \leq i_1 < \cdots < i_m \leq n$ of $W$, counted from the bottom of $W$, such that the $m \times m$ submatrix $W_m$ consisting of the rows $i_1, \ldots, i_m$ and columns $1, \ldots, m$ of $W$, has rank $m$.

Construct the $2^m$ elements of $V^m$ in lexicographic order: $V^m = \{0 = A_0, \ldots, A_{2^m-1}\}$. Precisely, $A_k = WK$, where $K$ is the column vector of the binary value of $k \in \{0, \ldots, 2^m - 1\}$, with the bits of $K$ in positions $i_1, \ldots, i_m$, the least significant bit in $i_1$; and 0's in all the other $n - m$ positions. See Table 1 for an example with $m = n$.

For given $d \in \mathbb{Z}^+$, scan $V^m$ from $A_0$ to $A_{2^m-1}$ to generate a subset $V' \subseteq V^m$ using the following greedy algorithm. Put $V' \leftarrow 0$. If $0 = A_{i_0}, \ldots, A_{i_j}$ have already been inserted into $V'$, insert $A_{i_{j+1}}$ if $i_{j+1} > i_j$ is the smallest integer such that $H(A_{i_\ell}, A_{i_{j+1}}) \geq d$ for $\ell \in \{0, \ldots, j\}$, where $H$ denotes Hamming distance.

TABLE 1. Generating a lexicode.

| $k$ | $V^m$ BIN | $V^m$ DEC | $V'$ |
|---|---|---|---|
| 0 | 0000 | 0 | ✓ |
| 1 | 0001 | 1 | |
| 2 | 0011 | 3 | ✓ |
| 3 | 0010 | 2 | |
| 4 | 0110 | 6 | ✓ |
| 5 | 0111 | 7 | |
| 6 | 0101 | 5 | ✓ |
| 7 | 0100 | 4 | |
| 8 | 1100 | 12 | ✓ |
| 9 | 1101 | 13 | |
| 10 | 1111 | 15 | ✓ |
| 11 | 1110 | 14 | |
| 12 | 1010 | 10 | ✓ |
| 13 | 1011 | 11 | |
| 14 | 1001 | 9 | ✓ |
| 15 | 1000 | 8 | |

The resulting $V'$ is the *lexicode* generated by $W$.

We remark that in Brualdi and Pless [1993], the term "lexicode" is reserved for the code generated when $W$ is the identity matrix, which is the case considered in Conway and Sloane [1986]; and "greedy codes" is used for the codes derived from any $W$ whose columns constitute a basis. Actually, in both of these papers no matrices are used, but the ordering is done in an equivalent manner. It seems natural, in the current context, to use matrices (see the proofs in the next section) and "lexicode" for the entire class of codes.

EXAMPLE 4. Let

$$W = \begin{pmatrix} \overset{4}{1} & \overset{3}{0} & \overset{2}{0} & \overset{1}{0} \\ 1 & 1 & 0 & 0 \\ 0 & 1 & 1 & 0 \\ 0 & 0 & 1 & 1 \end{pmatrix},$$

and $d = 2$, $m = n = 4$. We then get the ordered vector space depicted in Table 1. The vectors checked off in column $V'$ have been selected by our greedy algorithm, and constitute the lexicode.

Note that this lexicode is precisely the same code as that found in Example 1 by using a small groundgraph with $O(n^2)$ operations rather than $O(2^n)$ for the lexicode.

EXAMPLE 5. Let
$$W = \begin{pmatrix} 1 & 0 & 0 & 0 \\ 1 & 0 & 1 & 0 \\ 0 & 1 & 1 & 0 \\ 0 & 1 & 1 & 1 \end{pmatrix},$$
and $d = 2$. The reader should verify that the lexicode generated by $W$ is $(0, 7, 12, 11)$, in this order, which is identical to the code generated in Example 2.

EXAMPLE 6. Let
$$W = \begin{pmatrix} 1 & 0 & 0 & 0 & 0 \\ 1 & 1 & 0 & 0 & 0 \\ 0 & 1 & 1 & 0 & 0 \\ 0 & 0 & 1 & 1 & 0 \\ 0 & 0 & 0 & 1 & 1 \end{pmatrix},$$
and $d = 3$. The vector space now contains 32 entries, too large to list here. But the reader can verify that the lexicode generated by $W$ is precisely the same as that generated by the two polynomial methods of Example 3.

At the beginning of this section we defined $P$- and $N$-positions on the gamegraph $G$ of a game. We close the section by defining the *Sprague-Grundy* function $g = V \to \mathbb{Z}^0$ on $G(V, E)$ by $g(u) = \text{mex}\, g(F(u))$, where for any finite subset $S \subseteq \mathbb{Z}^0$,
$$\text{mex}(S) = \min(\mathbb{Z}^0 - S), \qquad g(S) = \{g(s) : s \in S\}.$$
It exists uniquely on any finite acyclic digraph. See, e.g., Berge [1985, 1989], Conway [1976] and Berlekamp, Conway and Guy [1982]. (Replace $g$ by $\gamma$ if $G$ has cycles or loops.)

**3. The Truth About Lexicodes.** With a lexicode in $V^m$, as defined in §2, associate a digraph $G = (V, E)$, called *lexigraph*, where $V$ is the set of all elements (vectors) of $V^m$; and $(A_k, A_j) \in E$ if and only if $j < k$ and $H(A_j, A_k) < d$, i.e., $w(A_j \oplus A_k) < d$, where for any vector $A$, $w(A) = weight$ of $A = H(A, 0) =$ number of 1-bits of $A$, where $\oplus$ denotes addition over $GF(2)$. If $(A_k, A_j) \in E$, then $A_j \in F(A_k)$ in the notation introduced at the beginning of §2. Note that $G$ is finite and acyclic. (For other possibilities of orienting the lexigraph, see the homework problem towards the end of this section.)

If we are so inclined, we can play a *lexigame* on $G$. Place a token on any vertex. A move consists of sliding the token from a vertex to a neighboring vertex along a directed edge. The player first unable to play loses and the opponent wins. The $P$-positions of the lexigame constitute the lexicode. (The lexigame is not overly interesting. The reason is that the lexigraph is "analogous" to a game-graph of a (more interesting) game played on a logarithmically smaller groundgraph. A game-graph of a game is not normally constructed but used instead for reasoning about the game. In fact, we do this in the proof of Theorem 3 below.)

We also point out that for lexicodes *per se* it suffices to consider the case $m = n$. It is only in §4, where we apply a greedy algorithm on anncodes, that the case $m < n$ will be important.

For any positive integer $s$, let $s^h$ denote the bit in the $h$-th digital position of the binary expansion of $s$, where $s^0$ denotes the least significant bit.

LEMMA 1. Let $a_1, a_2 \in \mathbb{Z}^0$, and let $b \in \{0, \ldots, a_1 \oplus a_2 - 1\}$. Then there is $i \in \{1, 2\}$ and $d \in \{0, \ldots, a_i - 1\}$ such that $b = a_j \oplus d$, $j \neq i$.

PROOF. Write $c = a_1 \oplus a_2$. Let $k = \max\{h : b^h \neq c^h\}$. Since $b < c$, we have $b^k = 0, c^k = 1$. Hence there exists $i \in \{1, 2\}$ such that $a_i^k = 1$. Letting $d = a_i \oplus b \oplus c = a_j \oplus b$, we have $d \in \{0, \ldots, a_i - 1\}$, since $b^h = c^h$ implies $d^h = a_i^h$ for $h > k$, and $d^k = 0$. ∎

For any $a \in \mathbb{Z}^0$, write $\phi(a) = \{0, \ldots, a - 1\}$.

COROLLARY 1. In the notation of Lemma 1, $0 \leq b < a_1 \oplus a_2$ implies $b \in a_1 \oplus \phi(a_2) \cup \phi(a_1) \oplus a_2$. ∎

By the closure of $V^m$, $\forall j \forall k \exists \ell$ such that $A_j \oplus A_k = A_\ell$.

LEMMA 2. We have $A_j \oplus A_k = A_{j \oplus k}$.

PROOF. As noted above, $A_j \oplus A_k = A_\ell$ for some $\ell$. Then $A_j = WJ$, $A_k = WK$, $A_\ell = WL$. Thus

$$WL = A_\ell = A_j \oplus A_k = W(J \oplus K).$$

This matrix equation implies $W_m L_m = W_m(J_m \oplus K_m)$, where $W_m$ was defined in §2, and any $m \times 1$ vector $H_m$ is obtained from the $n \times 1$ vector $H$ by retaining only the rows $i_1, \ldots, i_m$ of $H$, deleting the $n - m$ remaining rows, which contain only 0's for $L, J$ and $K$. Since $W_m$ is invertible, we thus get $L_m = J_m \oplus K_m$, so $\ell = j \oplus k$. ∎

The following is the main lemma of this section.

LEMMA 3. Let $A_j, A_k \in V^m$. Then for the lexigraph on $V^m$,

$$F(A_j \oplus A_k) \subseteq A_j \oplus F(A_k) \cup F(A_j) \oplus A_k \subseteq F(A_j \oplus A_k) \cup F^{-1}(A_j \oplus A_k).$$

PROOF. Let $A_\ell \in F(A_j \oplus A_k)$. By Lemma 2, $A_\ell \in F(A_{j \oplus k})$, so $w(A_\ell \oplus A_{j \oplus k}) = w(A_{j \oplus k \oplus \ell}) < d$ and $\ell < j \oplus k$. By Corollary 1, $\ell \in j \oplus \phi(k) \cup \phi(j) \oplus k$. Thus either there is $k' < k$ such that $\ell = j \oplus k'$, or there is $j' < j$ such that $\ell = j' \oplus k$. In the former case, $w(A_{j \oplus k \oplus \ell}) = w(A_{k \oplus k'}) < d$, so $A_\ell = A_j \oplus A_{k'} \in A_j \oplus F(A_k)$, and in the latter case we obtain, similarly, $A_\ell \in F(A_j) \oplus A_k$, establishing the left inclusion.

Let now $A_\ell \in A_j \oplus F(A_k) \cup F(A_j) \oplus A_k$. Then either $A_\ell = A_j \oplus A_{k'}$ for some $k' < k$ with $w(A_{k \oplus k'}) < d$, or $A_\ell = A_{j'} \oplus A_k$ for some $j' < j$ with $w(A_{j \oplus j'}) < d$. Without loss of generality, assume the former. Then $\ell = j \oplus k'$. Thus $w(A_{k \oplus k'}) = w(A_{j \oplus k \oplus \ell}) < d$. If $\ell < j \oplus k$, then $A_\ell \in F(A_j \oplus A_k)$, and if $\ell > j \oplus k$, then $A_j \oplus A_k \in F(A_\ell)$. ∎

We now show that the $g$-function is linear on the lexigraph $G$.

THEOREM 1. Let $G = (V, E)$ be a lexigraph. Then $g(u_1 \oplus u_2) = g(u_1) \oplus g(u_2)$ for all $u_1, u_2 \in V$.

PROOF. We write $(v_1, v_2) \in \mathcal{F}(u_1, u_2)$ if either $v_1 = u_1$ and $v_2 \in F(u_2)$; or else if $v_1 \in F(u_1)$ and $v_2 = u_2$, i.e.,

$$(v_1, v_2) \in \mathcal{F}(u_1, u_2) \quad \text{if} \quad (v_1, v_2) \in \big(u_1, F(u_2)\big) \cup \big(F(u_1), u_2\big). \tag{1}$$

Incidentally note that $\mathcal{F}$ is not a follower in $G$, but rather a follower in the sum game played on $G + G$. Let

$$K = \big\{(u_1, u_2) \in V \times V : g(u_1 \oplus u_2) \neq g(u_1) \oplus g(u_2)\big\},$$

$$k = \min_{(u_1, u_2) \in K} \big(g(u_1 \oplus u_2), g(u_1) \oplus g(u_2)\big).$$

If there is $(u_1, u_2) \in K$ such that $g(u_1 \oplus u_2) = k$, then $g(u_1) \oplus g(u_2) > k$. By Corollary 1 and the mex property of $g$, there is $(v_1, v_2) \in \mathcal{F}(u_1, u_2)$ such that $g(v_1) \oplus g(v_2) = k$. Now (1) implies

$$v_1 \oplus v_2 \in u_1 \oplus F(u_2) \cup F(u_1) \oplus u_2 \subseteq F(u_1 \oplus u_2) \cup F^{-1}(u_1 \oplus u_2),$$

where the inclusion follows from Lemma 3. Since $g(u_1 \oplus u_2) = k$, it follows that $g(v_1 \oplus v_2) > k$, so $(v_1, v_2) \in K$. Let

$$L = \big\{(u_1, u_2) \in K : g(u_1) \oplus g(u_2) = k\big\}.$$

We have just shown that $K \neq \emptyset$ implies $L \neq \emptyset$.

At this stage we note that the $\gamma$-function is a generalization of the $g$-function, so every $g$-function is also a $\gamma$-function. With the latter we can associate a monotonic counter function $c : V \to \mathbb{Z}^+$. We now pick $(u_1, u_2) \in L$ with $c(u_1) + c(u_2)$

– 9 –

minimum. For $(u_1, u_2) \in L$ we have $g(u_1 \oplus u_2) > k$. Then there is $v \in F(u_1 \oplus u_2)$ with $g(v) = k$. By the first inclusion of Lemma 3, there exists $(v_1, v_2) \in \mathcal{F}(u_1, u_2)$ such that $v = v_1 \oplus v_2$. So $g(v_1 \oplus v_2) = k$. Since $g(u_1) \oplus g(u_2) = k$, (1) implies $g(v_1) \oplus g(v_2) > k$, hence $(v_1, v_2) \in K$. As we saw earlier, this implies that there is $(w_1, w_2) \in \mathcal{F}(v_1, v_2)$ such that $(w_1, w_2) \in L$. Moreover, by property **B** of $\gamma$ (Definition 1 in Fraenkel and Yesha [1986]), we can select $(w_1, w_2)$ such that $c(w_1) + c(w_2) < c(u_1) + c(u_2)$, contradicing the minimality of $c(u_1) + c(u_2)$. Thus $L = K = \emptyset$. ∎

Let $V_i = \{u \in V : g(u) = i\}$ for $i \geq 0$. We now state the main result of this section.

THEOREM 2. Let $G = (V, E)$ be a lexigraph. Then $V_0 = V'$. Moreover, $V_0$ is a linear subspace of $V$. In fact, $g$ is a homomorphism from $V$ onto $GF(2)^t$ for some $t \in \mathbb{Z}^0$ with $t \leq \log_2 \left( \sum_{i=1}^{d-1} \binom{n}{i} \right)$, kernel $V_0$ and quotient space $V/V_0 = \{V_i : 0 \leq i < 2^t\}$, $\dim(V) = m + t$, where $m = \dim(V_0)$.

PROOF. By definition, $V$ is a vector space over $GF(2)$. Let $t$ be the smallest nonnegative integer such that $g(u) \leq 2^t - 1$ for all $u \in V$. Thus, if $t \geq 1$, there is some $v \in V$ such that $g(v) \geq 2^{t-1}$. Then the "1's complement" $2^t - 1 - g(v) < g(v)$. By the mex property of $g$, there exists $w \in F(v)$ such that $g(w) = 2^t - 1 - g(v)$. By Theorem 1, $g(v \oplus w) = g(v) \oplus g(w) = 2^t - 1$. Thus, again by the mex property of $g$, every value in $\{0, \ldots, 2^t - 1\}$ is assumed as a $g$-value by some $u \in V$. This last property holds trivially also for $t = 0$. Hence $g$ is onto. It is a homomorphism $V \to GF(2)^t$ by Theorem 1, and since $g(1u) = g(u) = 1g(u)$, $g(0u) = g(0 \cdots 0) = 0 = 0g(u)$.

By a standard result of linear algebra, $GF(2)^t \simeq V/V_0$, where $V_0$ is the kernel. Hence $V_0$ is a subspace of $V$. Clearly $V_0$ is also a graph-kernel of $G$. So is $V'$, which, by its definition, is both independent and dominating. Since any finite acyclic digraph has a unique kernel, $V_0 = V'$. Let $m = \dim(V_0)$. Then $\dim(V) = m + t$. The elements of $V/V_0$ are the cosets $V_i = w \oplus V_0$ for any $w \in V_i$ and every $i \in \{0, \ldots, 2^t - 1\}$. Finally, the outdegree of any vertex of $G$ is at most $\sum_{i=1}^{d-1} \binom{n}{i}$, so $2^t - 1 \leq \sum_{i=1}^{d-1} \binom{n}{i}$, which implies the bound on $t$. ∎

HOMEWORK 1. The lexigraph $G = (V, E)$ seems to exhibit a certain robustness, roughly speaking, with respect to $E$. That is, Theorem 2 seems to be invariant under certain edge deletions or reversions. In this direction, prove that Theorem 2 is still valid if $E$ is defined as follows: $(A_k, A_j) \in E$ if and only if $A_j < A_k$ (rather than $j < k$) and $H(A_j, A_k) < d$.

The proof of Theorem 2 is actually a much simplified version of the same result for annihilation games (Fraenkel and Yesha [1982]), where also the linearity of $\gamma$ was proved first. The simplification in the proof is no accident, since the lexigame played on the lexigraph (the groundgraph), can be considered to be an anngame with a single token. It's an acyclic groundgraph, for which the anngame theory is much simpler than for cyclic digraphs. We have given here a separate proof for lexicodes only because of the previous interest in them.

THEOREM 3. *The set of lexicodes is a subset of the set of anncodes.*

PROOF. Let $C$ be a lexicode with a given minimal distance. As we saw at the beginning of this section, $C$ is the set of the $P$-positions of the lexigame played on the lexigraph $G$. The lexigame is played on $G$ by sliding a token, and as such it is an annihilation game, whose $P$-positions constitute an anncode. Thus $C$ is an anncode. ∎

It might be of interest to explore the subset of anncodes generated when several tokens are distributed initially on a lexigraph, rather than only one.

Another question is under what conditions and for what finite fields $GF(p^a)$, $p$ prime, $a \in \mathbb{Z}^+$, are there "anncodes"? The key seems to be to generalize annihilation games as follows. On a given finite digraph $\Gamma$, place nonzero "particles" (elements of $GF(p^a)$), at most one particle per vertex. A move consists of selecting an occupied vertex and moving its particle to a neighboring vertex $v$ along a directed edge. If $v$ was occupied, then the "collision" generates a new particle, possibly 0 ("annihilation"), according to the addition table of $GF(p^a)$. The special case $a = 1$, when the particles are $0, \ldots, p - 1$, reduces to $p$-*annihilation*, i.e., the collision of particles $i$ and $j$ results in particle $k$, where $k \equiv i + j \pmod{p}$, $k < p$; and this special case becomes anngames for $p = 2$. Such "Elementary Particle Physics Games", whose $P$-positions are *collections* of linear codes, thus constitute a generalization of anngames. These games and their applications to coding seems to be an, as yet, unexplored area.

**4. Computing Anncodes.** In this section we give one particular example illustrating the computation of large anncodes. One can easily produce many other examples. The present example also shows how anncodes and lexicodes can be made to join forces.

We begin with a family $\Gamma_t$ of groundgraphs, which is a slightly simplified version of a family considered by Yesha [1978] for showing that the finite $\gamma$-values on an annihilation game played on a digraph without leaves can be arbitrarily large.

Let $t \in \mathbb{Z}^+$. We use the notation $J = J(t) = 2^{t-1}$, and define the digraph $\Gamma_t$ by

$$V(\Gamma_t) = \{x_i : 1 \leq i \leq J\} \cup \{y_i : 1 \leq i \leq J\}, \quad F(x_i) = y_i, \quad (i \in \{1, \ldots, J\}),$$
$$F(y_k) = \{y_i : 1 \leq i < k\} \cup \{x_j : 1 \leq j \leq J, j \neq k\} \quad (k \in \{1, \ldots, J\}).$$

See Fig. 4 for $\Gamma_3$.

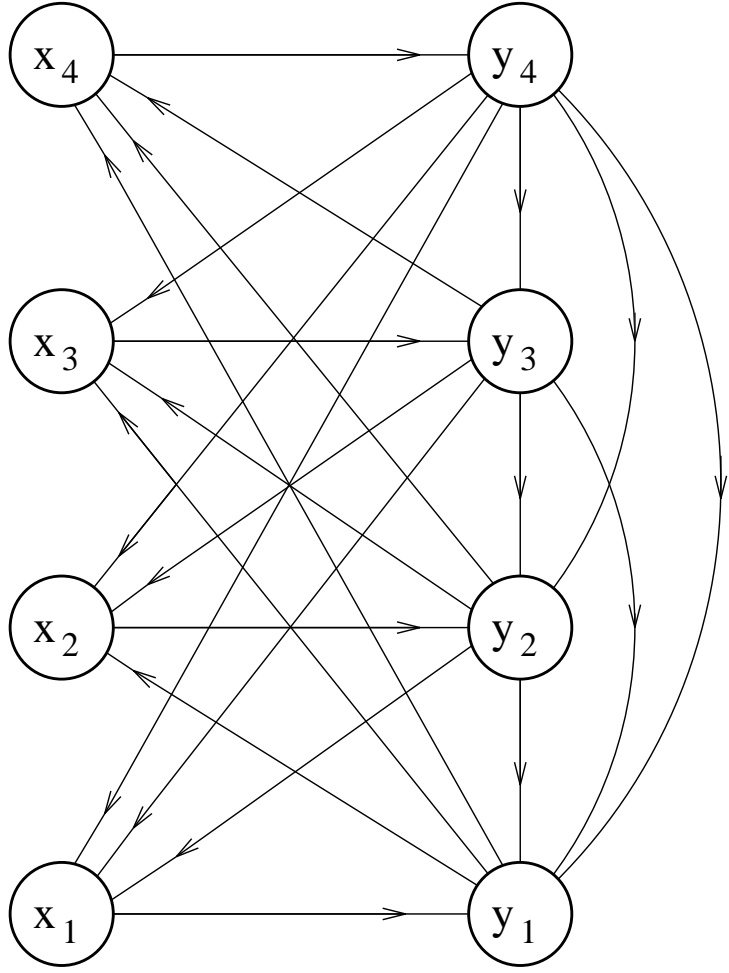

FIG. 4. The cyclic groundgraph $\Gamma_3$.

Since $\Gamma_t$ has no leaf, $\gamma(x_i) = \gamma(y_i) = \infty$ for all $i \in \{1, \ldots, J\}$. The following facts about the anngraph $G_t = (V, E)$ of $\Gamma_t$ are easy to establish, where $V^f = \{\boldsymbol{u} \in V : \gamma(\boldsymbol{u}) < \infty\}$.

(i)     $\gamma(x_i \oplus x_j) = 0$ for all $i \neq j$.

(ii)    $\gamma(x_i \oplus y_j) = j$ for all $i, j \in \{1, \ldots, J\}$.

(iii)   $\gamma(y_i \oplus y_j) = i \oplus j$ for all $i \neq j$.

(iv)   $\max_{\gamma(\boldsymbol{u}) < \infty} \gamma(\boldsymbol{u}) = \gamma(y_{J-1} \oplus y_J) = (J-1) \oplus J = 2^t - 1$.

(v)    $V^f = \{\boldsymbol{u} \in V : w(\boldsymbol{u}) \equiv 0 \pmod{2}\}$,   $|V^f| = 2^{2J-1}$,   $\dim(V^f) = 2J - 1$.

Thus, in the notation of Theorem 2, $m + t = 2J - 1$, hence

$$m = 2^t - t - 1.$$

– 12 –

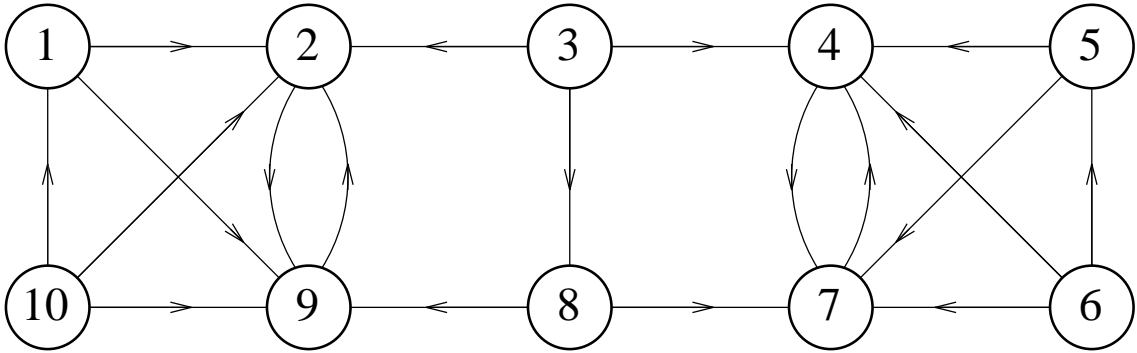

Fig. 5. The cyclic groundgraph $\Gamma'$.

For the family $\Gamma_t$ of groundgraphs, the $O(n^6)$ algorithm for computing $\gamma$ thus reduces to an $O(1)$ algorithm.

Now consider the groundgraph $\Gamma'$ depicted in Fig. 5. It is not hard to see that a basis for $V_0$ is given by the 4 vectors $(1,2,9,10)$, $(4,5,6,7)$, $(2,3,8,9)$, $(3,4,7,8)$. Each vector indicates the 4 vertices occupied by tokens.

We propose to play an annihilation game, say on $\Gamma = \Gamma_5 + \Gamma'$, which contains $32 + 10 = 42$ vertices. The vector space associated with the anngraph of $\Gamma$ contains $2^{42}$ elements, and to find a lexicode on $V^{42}$, for any given $d$, involves $2^{42}$ operations. On the other hand, we have for $\Gamma$, $|V_0| = m = 2^5 - 5 - 1 + 4 = 30$, so the anncode defined by $V_0$, for which $d = 2$, has $2^{30}$ elements. By the results of §2, we can compute a *lexi-anncode* for any $d > 2$, by applying the greedy algorithm to a lexicographic ordering of $V_0$, which can be obtained by using any basis of $V_0$. This computation involves only $2^{30}$ operations.

HOMEWORK 2. Carry out this computation, and find lexi-anncodes for several $d > 2$ on $\Gamma = \Gamma_5 + \Gamma'$.

We end with two further remarks.

(i) The Hamming distance between any two consecutive $P$-positions in an annihilation game is obviously $\leq 4$. Thus $d = 2$ for $\Gamma_t$ and $d = 4$ for $\Gamma'$. For finding codes with $d > 4$, it is thus natural to apply the greedy algorithm to a lexicographic ordering of $V_0$. Another method to produce anncodes with $d > 4$ is to encode each vertex of the groundgraph, i.e., each bit of the anngraph, by $k$ bits for some fixed $k \in \mathbb{Z}^+$. For example, in a lexigraph, each vertex is encoded by $n$ bits, and the distance between any two codewords is $\geq d$. In an Elementary Particle Physics game over $GF(p^a)$, it seems natural to encode each *particle* by $a$ digits. A third method for producing anncodes with $d > 4$ directly, seems to be to consider a generalization of anngames to the case where a move consists of sliding precisely $k$ (or $\leq k$) tokens, where $k$ is a fixed positive integer parameter — somewhat analogously to Moore's Nim (see e.g., Berlekamp, Conway and Guy

[1982, ch. 15]).

(ii) Note that $\bigcup_{i=0}^{2^k-1} V_i$ is a linear subspace of $V^f$ for every $k \in \{0,\ldots,t\}$. Any of these subspaces is thus also a linear code, in addition to $V_0$.

**Acknowledgment.** This paper is a direct result of Vera Pless' lecture at the Workshop on Combinatorial Games held at MSRI, Berkeley, CA, in July, 1994. I had intended to write it after I heard her lecture at the AMS Short Course on Combinatorial Games held in Columbus, OH, in the summer of 1990. But I put it off. This time, at $H^s$Lordship's banquet in Berkeley, Herb Wilf challenged me with his insights and comments, rekindling my interest in this topic. I also had a shorter conversation about it with John Conway and Bill Thurston at MSRI. Actually, I had originally discussed the possibility of using codes derived from annihilation games with C.L. Liu and E.M. Reingold in 1978 or 1979, and again briefly with Richard Guy in the 80's and with Ya'acov Yesha in the early 90's.